\documentclass[12pt]{amsart}
\input epsf
\usepackage{amsmath}
\usepackage{amssymb}
\usepackage{amscd}
\usepackage{amsthm}

\theoremstyle{plain}
\newtheorem{theorem}[equation]{Theorem}

\theoremstyle{remark}

\theoremstyle{definition}

\newcommand{\ra}{\rightarrow}

\newcommand{\R}{\mathbb R}

\newcommand{\K}{{\mathcal K}}

\newcommand{\M}{{\mathcal M}}

\newcommand{\genus}{\operatorname{genus}}

\newcommand{\Int}{\operatorname{Int}}

\def\D{\partial}
\newcommand{\al}{\alpha}
\def\de{\delta}

\def\eps{\epsilon}
\def\ga{\gamma}

\def\ra{\rightarrow}

\def\defeq{:=}

\begin{document}

\title{Singularity structure in mean curvature flow of mean convex
sets}

\author{Tobias H. Colding}
\thanks{THC was supported by NSF grant
DMS 0104453.}

\email{colding@cims.nyu.edu}

\address{Courant Institute of Mathematical Sciences\\
        251 Mercer Street\\
        New York, NY 10012}

\author{Bruce Kleiner}
\thanks{BK was supported by NSF grant
DMS-0204506.}

\email{bkleiner@umich.edu}

\address{Department of Mathematics\\
        University of Michigan\\
        2072 East Hall, 525 E Univ. Ave.\\
        Ann Arbor, MI 48109-1109}
\date{October 15, 2003}
\maketitle

\begin{abstract}
In this note we announce results on the mean curvature flow of
mean convex sets in 3-dimensions.  Loosely speaking, 
our results justify the naive picture of
mean curvature flow where the only singularities are neck
pinches, and components which collapse to asymptotically round
spheres.  
\end{abstract}

In this note we announce results on the mean curvature flow of
mean convex sets; all the statements below have natural generalizations
to the setting of Riemannian $3$-manifolds, but for the sake of simplicity
we will primarily discuss subsets of $\R^3$ here. Loosely speaking, 
our results justify the naive picture of
mean curvature flow where the only singularities are neck
pinches, and components which collapse to asymptotically round
spheres.  
Recall that a one-parameter family of smooth hypersurfaces $\{ M_t \} \subset
\R^{n+1}$ {\it flows by mean curvature} if
\begin{equation}\label{e:meanflow}
    z_t = {\bf{H}} (z) = \Delta_{M_t} z \, ,
\end{equation}
where $z=(z_1,\ldots,z_{n+1})$ are coordinates on $\R^{n+1}$ and 
${\bf{H}} = - H {\bf{n}}$
is the mean curvature vector.
The papers 
\cite{evansspruck} and \cite{chengigagoto} defined a level set flow for
any closed subset $K$ of $\R^n$.  This is a $1$-parameter family
of closed sets $K_t\subset \R^n$ with $K_0=K$ (when $K$ is
a domain bounded by a smooth compact hypersurface then the evolution
of $\D K$ for a short time interval coincides with the classical
mean curvature evolution).  Following \cite{white1},
we say that a compact subset $K\subset \R^n$ is
{\em mean convex} if $K_t\subset\Int(K)$ for all $t>0$.
In this case there is also an associated Brakke flow $\M:t\mapsto M_t$
 of rectifiable varifolds \cite{brakke,ilmanen,white1}, and
the pair $(\M,\K)$, where 
$$\K\defeq \bigcup_{t\geq 0}\,K_t\times\{t\}\subset 
\R^n\times\R$$
 is called a {\em mean-convex flow}, \cite{white2}.
The fundamental papers \cite{white1,white2} developed a far-reaching
partial regularity theory for mean curvature flow of mean convex
subsets of $\R^n$.  Our results build on \cite{white1,white2}, giving
finer understanding of the singularities 
in the $3$-dimensional case.  Recall that the main result of \cite{white1} 
asserts that the space time singular set of the region swept out by a 
mean--convex set in $\R^{n+1}$ has parabolic Hausdorff dimension at most 
$(n-1)$, and  \cite{white2} proved a structure theorem for 
blow--ups of mean--convex flows; cf. also  
\cite{huiskensinestrari1,huiskensinestrari0}.  
We expect that the more refined description
 of singularities given here 
will open the way for applications of mean convex flow to geometric
and/or topological problems involving mean convex surfaces.

When $(\M,\K)$ is a mean convex flow in $\R^3$,
then for almost every time $t$ the time
slice $K_t$ is a domain with
 smooth boundary, \cite[Corollary to Theorem 1.1]{white1}.  Our first result 
shows that the high curvature
portion of such smooth time slices has standard local geometry:  
\begin{theorem}
For all $\eps > 0$ there is a number $h_0=h_0(\eps)$ with the
following property.  If $(\M,\K)$ is a mean convex flow in $\R^3$
and  $K_t$ is a regular time slice of $\K$ for some $t>0$, 
then there is a decomposition
$K_t=G_t\cup B_t$, such that

\begin{itemize}
\item For all $x\in G_t$,
and after rescaling by the factor $\frac{h_0}{d(x,\D K)}$ 
the pointed subset $(K_t,x)$ is  $\eps$-close
to some pointed half-space $(P,p)$ in the pointed
$C^{\frac{1}{\eps}}$-topology.

\medskip
\item Each component of $B_t$ is diffeomorphic to the $3$-ball or a solid torus, and for all
$x\in \D K_t\cap B_t$, the pointed subset $(K_t,x)$ becomes, after rescaling
by the factor $H(x)$, $\eps$-close to a pointed convex model subset $(V,v)$
in the pointed $C^{\frac{1}{\eps}}$-topology.  Here $V\subset \R^3$ is a 
convex set whose
tangent cone at infinity is either a point, a line, or a ray, and $V$
looks like a round cylinder near infinity, in the following sense: 
for every $\de>0$ there is a compact set $K\subset V$, such that
for every  $v'\in V$ lying 
 outside $K$, if we rescale $V$ by $H(v')$,
the resulting pointed subset $(V,v')$ is $\de$-close to a round
cylinder in the pointed $C^{\frac{1}{\de}}$-topology.
\end{itemize}
\end{theorem}
\noindent
 Note that the bounds on the geometry deteriorate
as one approaches $\D K$; this is by necessity since no regularity
condition has been imposed on $K$.  If $K$ happens to be smooth,
then standard estimates for smooth mean curvature flow 
control the geometry of $K_t$ when $t\lesssim \sqrt r$, where
$r$ is the normal injectivity radius of $\D K$.
Theorem 2 may be compared with the recent work of 
Huisken-Sinestrari \cite{huiskensinestrari},
where a similar geometric description was obtained for mean curvature 
flow of smooth hypersurfaces in $\R^n$
where the sum of the first two 
principal curvatures is positive.  
The results in 
\cite[sections 11, 12]{perelman1} are also in a similar spirit.  Note that
their results only apply to the evolution prior to the formation of the first
singularity, whereas our results, like those in \cite{white1,white2},
apply even after the formation of
a singularity.  (In fact,
the methods yield a decomposition of arbitrary time slices, which we
omit for the sake of simplicity.)

It follows from the strong maximum principle and compactness that the sets $\D K_t$ for $t\geq 0$ are disjoint, and define a ``singular
foliation'' of the original set $K$.   Our next theorem proves H\"older regularity of the singular set of the foliation $\partial K_t$.  

\begin{theorem}

The foliation defined by the sets $\D K_t$
is smooth on the complement of a closed subset $S\subset K$
which satisfies the following Reifenberg-type condition: for all $\eps>0$  
there is an $r_0=r_0(\eps)$ such that if $r<r_0$
and $x\in S$, then there is a line 
$A\subset \R^3$ such that $S\cap B(x,r)$ is contained in the
tubular neighborhood $N_{\eps r}(A)$.  In particular, $S$ lies
in a $1$-dimensional topological submanifold $\ga\subset K$ which admits
a $C^\al$-biH\"older parametrization for all $\al < 1$.  Furthermore, the mean curvature
defines a  proper function on $K\setminus S$.
\end{theorem}

After passing through a singularity 
the topological type of a surface 
flowing by mean curvature can change.  In \cite{white3} White proved some 
results comparing the homology of the surface before and after such a 
singularity.  Our next theorem shows that the region between two regular 
time slices is obtained
from the earlier time slice by attaching $2$ and $3$--handles.  
Recall that attaching a $k$--handle to the boundary of an $n$-manifold $N$
is essentially just the process of attaching a fattened-up 
$k$--disk to $\D N$ along the 
$(k-1)$--sphere, i.e. one glues $D^k\times D^{n-k}$ to $\D N$ 
along $\D D^k\times D^{n-k}$.

\begin{theorem}
If $0\leq t < t'$ and $K_t, K_{t'}$ are regular time slices, then 
$K_t\setminus \Int(K_{t'})$ is a compact $3$-manifold with boundary
which may be obtained from $\D K_t$ by attaching $k$-handles for 
$k=2,3$.  
\end{theorem}

Our final theorem deals with mean convex flow
in a general $3$-manifold, where the 
flow may converge as time tends to infinity to a  set $K_\infty$ with 
nonempty interior.

\begin{theorem}
Let $M$ be a compact Riemannian $3$-manifold, and $K\subset M$
a mean convex subset with smooth boundary.  Then as $t\ra\infty$,
the intersection of the sets $K_t$ converges to a (possibly empty)
domain $K_\infty\subset \Int(K)\subset M$, where each boundary
component of $K_\infty$ is a smooth, weakly stable minimal surface,
and $\genus(\D K_\infty)\leq \genus(\D K)$.   Furthermore, any 
compact minimal surface in $K\setminus \Int(K_\infty)$
is contained in $\D K_\infty$; in particular $\D K_\infty$
is homologically minimizing in the domain $K\setminus\Int(K_\infty)$.
\end{theorem}

\bibliography{refs}

\providecommand{\bysame}{\leavevmode\hbox to3em{\hrulefill}\thinspace}
\providecommand{\MR}{\relax\ifhmode\unskip\space\fi MR }
\providecommand{\MRhref}[2]{%
  \href{http://www.ams.org/mathscinet-getitem?mr=#1}{#2}
}
\providecommand{\href}[2]{#2}
\begin{thebibliography}{CGG91}

\bibitem[Bra78]{brakke}
K.~Brakke, \emph{The motion of a surface by its mean curvature}, Mathematical
  Notes, vol.~20, Princeton University Press, Princeton, N.J., 1978.

\bibitem[CGG91]{chengigagoto}
Y.~G. Chen, Y.~Giga, and S.~Goto, \emph{Uniqueness and existence of viscosity
  solutions of generalized mean curvature flow equations}, J. Differential
  Geom. \textbf{33} (1991), no.~3, 749--786.

\bibitem[ES91]{evansspruck}
L.~C. Evans and J.~Spruck, \emph{Motion of level sets by mean curvature. {I}},
  J. Differential Geom. \textbf{33} (1991), no.~3, 635--681.

\bibitem[HS]{huiskensinestrari}
G.~Huisken and C.~Sinestrari, in preparation.

\bibitem[HS99a]{huiskensinestrari0}
\bysame, \emph{Convexity estimates for mean curvature flow and singularities of
  mean convex surfaces}, Acta Math. \textbf{183} (1999), no.~1, 45--70.

\bibitem[HS99b]{huiskensinestrari1}
\bysame, \emph{Mean curvature flow singularities for mean convex surfaces},
  Calc. Var. Partial Differential Equations \textbf{8} (1999), no.~1, 1--14.

\bibitem[Ilm94]{ilmanen}
T.~Ilmanen, \emph{Elliptic regularization and partial regularity for motion by
  mean curvature}, Mem. Amer. Math. Soc. \textbf{108} (1994), no.~520, x+90.

\bibitem[Per02]{perelman1}
G.~Perelman, \emph{The entropy formula for the {R}icci flow and its geometric
  applications}, math.DG/0211159.

\bibitem[Whi95]{white3}
B.~White, \emph{The topology of hypersurfaces moving by mean curvature}, Comm.
  Anal. Geom. \textbf{3} (1995), no.~1-2, 317--333.

\bibitem[Whi00]{white1}
\bysame, \emph{The size of the singular set in mean curvature flow of
  mean-convex sets}, J. Amer. Math. Soc. \textbf{13} (2000), no.~3, 665--695.

\bibitem[Whi03]{white2}
\bysame, \emph{The nature of singularities in mean curvature flow of
  mean-convex sets}, J. Amer. Math. Soc. \textbf{16} (2003), no.~1, 123--138.

\end{thebibliography}
\bibliographystyle{amsalpha}

\end{document}